\documentclass[10pt,letterpaper,dvips]{article}
 \usepackage{amssymb,latexsym,amsmath,amsthm}
 
\usepackage{fullpage}
\usepackage{url}
\usepackage{xypic}
\usepackage[all]{xy}

\usepackage[numbers,sort&compress]{natbib}

\usepackage{wasysym}

\newtheorem{thm}{Theorem}
\newtheorem{lemma}{Lemma}

\newtheorem{cor}{Corollary}

\newtheorem{rem*}{Remark}

\newtheorem{conj}{Conjecture}

\begin{document}
\date{}

\title{Liouville theorems for the polyharmonic H\'{e}non-Lane-Emden system}
\maketitle

\begin{center}
\author{{Mostafa Fazly}
}
\\
{\it\small Department of Mathematical and Statistical Sciences, University of Alberta}\\
{\it\small Edmonton, Alberta, Canada T6G 2G1}\\
{\it\small e-mail: fazly@ualberta.ca}\vspace{1mm}
\end{center}

\vspace{3mm}

\begin{abstract} 
We study Liouville theorems for the following polyharmonic H\'{e}non-Lane-Emden system 
\begin{eqnarray*}
 \left\{ \begin{array}{lcl}
\hfill (-\Delta)^m  u&=& |x|^{a}v^p   \ \ \text{in}\ \ \mathbb{R}^n,\\   
\hfill (-\Delta)^m  v&=& |x|^{b}u^q   \ \ \text{in}\ \ \mathbb{R}^n,
\end{array}\right.
  \end{eqnarray*}
when $m,p,q \ge 1,$ $pq\neq1$, $a,b\ge0$.  The main conjecture states that $(u,v)=(0,0)$ is the unique nonnegative solution of this system whenever  $(p,q)$ is  {\it under} the critical Sobolev hyperbola, i.e.
$ \frac{n+a}{p+1}+\frac{n+b}{q+1}>{n-2m}$.    We show that this is indeed the case  in dimension $n=2m+1$ for bounded solutions.  In particular, when $a=b$ and $p=q$, this means that $u=0$ is the only nonnegative bounded  solution of the polyharmonic H\'{e}non equation 
 \begin{equation*}
 (-\Delta)^m  u= |x|^{a}u^p   \ \ \text{in}\ \ \mathbb{R}^{n}
  \end{equation*}
  in dimension $n=2m+1$ provided $p$ is the  subcritical Sobolev exponent, i.e., $1<p<{1+4m+2a}$.   Moreover, we show that the conjecture holds for radial solutions in any dimensions. It seems the power weight functions $|x|^a$ and $|x|^b$ make the problem dramatically more challenging when dealing with nonradial solutions.

\end{abstract}


\section {Introduction and main results}
We examine the following weighted system known as the polyharmonic H\'{e}non-Lane-Emden system 
\begin{eqnarray}
\label{main}
 \left\{ \begin{array}{lcl}
\hfill (-\Delta)^m u&=& |x|^{a}v^p   \ \ \text{in}\ \ \mathbb{R}^n,\\   
\hfill (-\Delta)^m v&=& |x|^{b}u^q   \ \ \text{in}\ \ \mathbb{R}^n,
\end{array}\right.
  \end{eqnarray}
where $m,p,q \ge 1$, $pq\neq1$, $a,b\ge0$.  This is the statement of the  {\it H\'{e}non-Lane-Emden conjecture} for polyharmoic system \eqref{main}. 
\begin{conj}\label{conj} Let $(u,v)$ be a nonnegative solution of system (\ref{main}).  Suppose $(p,q)$ is {\it under} the critical hyperbola, i.e.,
  \begin{equation}\label{curve}
   \frac{n+a}{p+1}+\frac{n+b}{q+1}>{n-2m}.
  \end{equation}
Then $u=v=0$. 
\end{conj}
Note that it is very straightforward to give a positive answer to Conjecture \ref{conj} in all dimensions $1\le n\le 2m$. This is in fact a quick consequence of $L^1$ estimates given in Lemma \ref{weightedest}.  Therefore, in this paper, we focus on higher dimensions $n\ge 2m+1$. 

Liouville theorems for system (\ref{main}) are very widely studied for the past few decades. In what follows we briefly review some of the related known results.   We divide the introduction into two different cases. We first consider the case  $a=b=0$ and then the other case that is when one of the parameters $a$ or $b$ is not zero.  In this paper, we mainly focus on system (\ref{main}) whenever there are weight functions $|x|^a$ and $|x|^b$. As a matter of fact, the weight functions make the problem much more challenging and as a general statement, some standard techniques such as moving plane methods and certain Sobolev embeddings cannot be applied anymore. 

\subsection{The case $a=b=0$}
 System (\ref{main}) when $a=b=0$ is well studied and there are many interesting results on classifying the solutions of this system for various $p$ and $q$.    
 
 We begin by the scalar case that is when $p=q>1$.  For the Lane-Emden equation (i.e., when $m=1$, $p=q>1$ and $a=b=0$) a celebrated theorem by Gidas-Spruck  \cite{gs,gs2}  states that there is no positive solution for the Lane-Emden equation 
\begin{equation*}
  -\Delta u= u^p\, \ \ \text{in} \ \ \mathbb{R}^n
  \end{equation*} 
   whenever $1<p< \frac{n+2}{n-2}$ for $n\ge 3$.   This Liouville theorem is  optimal as shown by Gidas, Ni and Nirenberg in \cite{gnn} under the assumption that $u=O(|x|^{2-n})$, and by Caffarelli, Gidas and Spruck in \cite{cgs} without  the growth assumption. See also Chen and Li \cite{cl}  for an easier proof based on the  moving planes method.  In the case of the fourth order Lane-Emden equation (i.e., when $m=2$, $p=q>1$ and $a=b=0$) and the polyharmonic Lane-Emden equation (i.e., when $m\ge 1$, $p=q>1$ and $a=b=0$) 
   \begin{equation*}
 ( -\Delta)^m u= u^p\, \ \ \text{in} \ \ \mathbb{R}^n
  \end{equation*} 
   similar Liouville theorems are proved by Lin \cite{lin} and Wei and Xu in \cite{wx} for the subcritical Sobolev exponent that is $1<p< \frac{n+2m}{n-2m}, n>2m$.  Note that this exponent appears in the Sobolev embedding $W^{m,2}\hookrightarrow L^p$. 

Now we focus on the case that the parameters $p$ and $q$ are not necessarily equal. Therefore, we are dealing with a system of equations. This case is much less understood than the scalar case.   For the Lane-Emden system (i.e., when $m=1$, $a=b=0$ and $p,q \ge1$ when $pq\neq1$)
\begin{eqnarray*}
 \left\{ \begin{array}{lcl}
\hfill -\Delta u&=& v^p   \ \ \text{in}\ \ \mathbb{R}^n,\\   
\hfill -\Delta v&=& u^q   \ \ \text{in}\ \ \mathbb{R}^n,
\end{array}\right.
  \end{eqnarray*}
 Conjecture \ref{conj} is known as the Lane-Emden conjecture and the curve $ \frac{1}{p+1}+\frac{1}{q+1}=\frac{n-2}{n}$ is the {\it critical Sobolev hyperbola}.    Proving such a nonexistence result seems to be very challenging problem.  However, there are many interesting papers that cover certain dimensions.     
  The case of radial solutions was solved by Mitidieri \cite{m} in any dimension, and both Mitidieri \cite{m} and Serrin-Zou \cite{sz98} constructed positive radial solutions {\it on} and {\it above} the critical hyperbola, i.e. $\frac{1}{p+1}+\frac{1}{q+1}\le\frac{n-2}{n}$, which means that the nonexistence theorem is optimal for radial solutions. 
  For nonradial solutions of the Lane-Emden system, certain Liouville theorems are proved for various parameters $p$ and $q$ by Souto in \cite{s95},  Mitidieri in \cite{m} and Serrin-Zou in \cite{sz}, that in a particular case they give a positive answer to the Lane-Emden conjecture in dimensions $1\le n\le 2$. In  dimension $n=3$,  Serrin-Zou \cite{ sz} gave a proof for the nonexistence of polynomially bounded solutions, an assumption that was relaxed later by Pol\'{a}\v{c}ik, Quittner and Souplet \cite{pqs}. See also \cite{bm,fd}.   More recently, Souplet \cite{so} completely settled the conjecture in dimension $n=4$, while providing in dimensions $n\ge 5$, a more restrictive new region for the exponents $(p,q)$ that insures nonexistence.  The Lane-Emden conjecture is an open problem for dimensions $n\ge 5$. For the polyharmonic case $m\ge 1$, it is known that  $ \frac{1}{p+1}+\frac{1}{q+1}=\frac{n-2m}{n}$ is the {\it critical Sobolev hyperbola}.   Conjecture \ref{conj} for the case of radial solutions was solved by  Liu et. al.  in \cite{lgz} in any dimensions and as far as we know only some partial results are given for the nonradial solutions in \cite{lgz,y,lei}.   Note that  Caristi,  DÕAmbrosio and Mitidieri  in  \cite{cam}  have proved Liouville theorems for supersolutions of system (\ref{main})  and also  they have explored the  connection between (\ref{main}) and the Hardy-Littlewood-Sobolev systems (HLS). 
 
\subsection{The case $a\neq0$ and or $b\neq0$}

  The power weight function has been of interest in this context and it was introduced by M. H\'{e}non \cite{h} in equation 
\begin{equation}\label{henon}
\left\{
                      \begin{array}{ll}
                       -\Delta u=|x|^a u^{p}  & \hbox{for  $|x|<1$,} \\
                       u=0 & \hbox{on $|x|=1$,} \\
                           \end{array}
                    \right.
                    \end{equation}  
to model and study spherically symmetric clusters of stars. This equation is now known as the H\'{e}non equation for $a>0$ and the H\'{e}non-Hardy equation for $a<0$.  Ten years later, Ni in \cite{n} explored properties of positive radial solutions of the H\'{e}non equation on the unit ball and observed the fact that the power profile $|x|^a$ enlarges considerably the range of solvability beyond the classical critical threshold, i.e., $p<2^*-1=\frac{n+2}{n-2}$ to $p<2^*-1+\frac{2a}{n-2}=\frac{n+2+2a}{n-2}$ where $2^*$ is the critical Sobolev exponent for the Sobolev embedding $W^{1,2}\hookrightarrow L^p$.   On the other hand, as it is shown by Smets, Su and  Willem in \cite{ssw} and references therein, equation (\ref{henon}) also admits nonradial solutions for $p<2^*-1$.  The existence of nonradial solutions for the full range $p<\frac{n+2+2a}{n-2}$ is still an open problem.     Note that since the function $|x|\to |x|^a$ is increasing, the classical moving planes arguments given by Gidas, Ni and Nirenberg   in \cite{gnn} cannot be applied to prove the radial symmetry of the solutions of (\ref{henon}). Therefore, the existence of nonradial solutions for this equation is natural and it is studied in many interesting papers. 

Regarding Liouville theorems, Phan and Souplet \cite{ps} for the first time attacked the problem and showed among other results that conjecture for the scalar case, that is $m=1$ and $a=b$ and $p=q>1$, holds for bounded nonnegative solutions in dimension $n=3$.   Here is the result,

  \medskip
\noindent{\bf Theorem A. (Phan-Souplet \cite{ps})}\quad
{\it Let $n=3$, $m=1$, $a=b >-2$ and $p=q>1$. Assume $(p,q)$ satisfies (\ref{curve}) that is $1<p<5+2a$, then there is no positive bounded solution for the H\'{e}non equation, i.e., 
\begin{equation}\label{mainsingle}
  -\Delta u=|x|^a u^p\, \ \ \text{in} \ \ \mathbb{R}^n.
  \end{equation} 
}
For the case of systems (\ref{main}) when $m=1$, the author with Ghoussoub in \cite{fg} have also proved the conjecture in dimension three for bounded solutions. See also \cite{phan} for dimensions three and some partial results for dimension four.   In this note, we shall first extend the result of Fazly-Ghoussoub \cite{fg} and  Phan-Souplet \cite{ps} to the full polyharmonic H\'{e}non-Lane-Emden system by showing the following.

\begin{thm}\label{result}
Conjecture \ref{conj} holds in dimension $n=2m+1$ for nonnegative bounded  solutions of (\ref{main}).  
\end{thm}

\begin{thm}\label{resultrad} Conjecture \ref{conj} holds in all dimensions for nonnegative radial solutions of (\ref{main}).  

\end{thm}

For the special case $p=q$ and $a=b$, we have  the following weighted equation known as the polyharmonic H\'{enon} equation 
\begin{equation} \label{maineq}
 (-\Delta)^m u = |x|^{a}u^p   \ \ \text{in}\ \ \mathbb{R}^n,
  \end{equation}
where $p >1$ and $a\ge0$. Note that {\it under} the critical hyperbola (\ref{curve}) turns into the following subcritical Sobolev exponent 
 \begin{equation}\label{curveeq}
  1<p<\frac{n+2m+2a}{n-2m} \ \ \text{where} \ \ n>2m
  \end{equation}
 As a consequence of Theorem \ref{result}, $u=0$ is the unique nonnegative bounded solution of (\ref{maineq})  in dimension $n=2m+1$ provided (\ref{curveeq}) holds that is $1<p<1+4m+2a$. Also, Theorem \ref{resultrad} implies that $u=0$ is the unique nonnegative radial solution of (\ref{maineq}) provided (\ref{curveeq}) holds in all dimensions. Let us mention that very recently, Cowan in \cite{cow} following ideas developed in \cite{ps,fg} considered the fourth order H\'{e}non equation, that is (\ref{maineq}) when $m=2$, and proved that in dimension five there is no bounded positive solution for (\ref{maineq}) provided $1<p<9+2a$.

Here is the organization of the paper.   In Section \ref{nonrad}, we prove Theorem \ref{result} via applying various methods developed  in the theory of elliptic regularity. In Section \ref{rad}, we prove Theorem \ref{resultrad} via certain ODE arguments. Our methods of proof are strongly motivated by the ideas developed by Souplet in \cite{so},  Phan-Souplet in \cite{ps}, Mitidieri et. al.  in \cite{mi,m,mp},  Wei-Xu in \cite{wx}, Fazly-Ghoussoub in \cite{fg} and references therein.

\section{Liouville theorems for nonradial solutions via elliptic estimates}\label{nonrad}

We start with the following standard $L^1(B_R)$ estimate on the right hand side of system (\ref{main}). Similar estimates for the second order case are given in  \cite{mp,qs,sz}. 

\begin{lemma} \label{weightedest}
For any nonnegative entire solution $(u,v)$ of (\ref{main}) and $R>1$, there holds  
\begin{eqnarray} \label{weightedest1}
\int_{B_R}{       |x|^{a}   v^p} &\le&  C \ R^{n -2m- \frac{(b+2m)p+(a+2m)}{pq-1} },
\\ \label{weightedest2}
\int_{B_R}{    |x|^{b}   u^q}  &\le& C\  R^{n-2m - \frac{(a+2m)q+(b+2m)}{pq-1} },
 \end{eqnarray}
where the positive constant $C=C(n,m,t,s,a,b,p,q)$ does not depend on $R$.
\end{lemma}

Note that (\ref{weightedest1}) and (\ref{weightedest2}) imply that if  $(u,v)$ is a nonnegative entire solution of (\ref{main}) then  $u=v=0$ in dimensions 
\begin{eqnarray*}
n-2m< \max\left\{  \frac{(b+2m)p+(a+2m)}{pq-1},  \frac{(a+2m)q+(b+2m)}{pq-1}\right\}.
 \end{eqnarray*}
 In particular this proves Conjecture \ref{conj} in dimensions $1\le n\le 2m$. However this does not cover {\it under} the critical hyperbola mentioned in Conjecture \ref{conj} for $n\ge 2m+1$. In other words, (\ref{curve}) is equivalent to the following when $pq>1$
\begin{eqnarray*}
n-2m<  \frac{(b+2m)(p+1)+(a+2m)(q+1)}{pq-1}.
 \end{eqnarray*}
 \textbf{Proof:} Fix the following standard test function $\phi_R\in C^\infty_c(\mathbb{R}^n)$ when $0\le\phi_R\le1$ and 
 $$\phi_R(x)=\left\{
                      \begin{array}{ll}
                        1, & \hbox{if $|x|<R$;} \\
                        0, & \hbox{if $|x|>2R$;} 
                                                                       \end{array}
                    \right.$$
where $|| D_x^i\phi_R||_{L^{\infty}(\mathbf R^n)} \le \frac{C}{R^i}$ for $i=1,\cdots,2m$. For any $t\ge 2m$, we have  
$$|\Delta^m \phi^t_R(x)|\le C \left\{
                      \begin{array}{ll}
                        0, & \hbox{if $|x|<R$ or $|x|>2R$;} \\
                         R^{-2m} \phi^{t-2m}_R , & \hbox{if $R<|x|<2R$;} 
                                                                       \end{array}
                    \right.$$
Now test the first equation of (\ref{main}) by $\phi^t_R$ and integrate to get 
\begin{eqnarray*}
 \int_{B_{2R}}     |x|^{a}     v^p \phi^t_R 
  \le C  R^{-2m}\int_{B_{2R}\setminus B_{R}} u\phi^{t-2m}_R.
\end{eqnarray*}
Applying the H\"{o}lder's inequality we get 
\begin{eqnarray*}
\int_{B_{2R}}      |x|^{a}     v^p  \phi^t_R
& \le  & C\  R^{  (n-\frac{b}{q}q')\frac{1}{q'}  -2 m}  \left(     \int_{B_{2R}\setminus B_{R}}       |x|^{b}     u^q \phi^{(t-2m)q}_R   \right)^{1/q}.
\end{eqnarray*}
By a similar calculation for $s\ge 2m$, we obtain
\begin{eqnarray*}
 \int_{B_{2R}}    |x|^{b}      u^q  \phi^s_R & \le &C \ R^{  (n-\frac{a}{p}p')\frac{1}{p'}  -2m }     \left(\int_{B_{2R}\setminus B_{R}}    |x|^{a} v^p \phi^{(s-2m)p}_R   \right)^\frac{1}{p} ,
 \end{eqnarray*}
where $\frac{1}{p}+\frac{1}{p'}=1$. Since $pq>1$, for large enough $s$ we have $2m+\frac{s}{q}<(s-2m)p$. So, we can choose $t$ such that $2m+\frac{s}{q}\le   t  \le   (s-2m)p$ which means that $t\le (s-2m)p$ and $s\le (t-2m)q$.  Therefore, $\phi^{(t-2m)q}_R \le \phi^{s}_R$  and  $ \phi^{(s-2m)p}_R \le  \phi^{t}_R$. Now, by collecting  the above inequalities  for $pq> 1$ we get
 
 \begin{eqnarray}
\nonumber   \left(\int_{B_{2R}}{         |x|^{a}      v^p   \phi^t_R} \right)^{pq}&\le & C\    R^{ [ (n-\frac{b}{q}q')\frac{1}{q'}  -2 m] pq}  \left(  \int_{B_{2R}\setminus B_{R}}  {    |x|^{b}  u^q \phi^s_R}  \right)^p 
\\ \label{L1final1}
&\le & C\ R^{   (n-2m)(pq-1) -[(b+2m)p+(a+2m)] }  \int_{B_{2R}\setminus B_{R}}     |x|^{a}      v^p   \phi^t_R,
 \end{eqnarray}
 and
 \begin{eqnarray}
  \nonumber       \left(\int_{B_{2R}} {    |x|^{b}  u^q     \phi^s_R } \right)^{pq} &\le & C\    R^{ [  (n-\frac{a}{p}p')\frac{1}{p'}  -2m  ] pq}   \left( \int_{B_{2R}\setminus B_{R}}      {         |x|^{a}      v^p   \phi^t_R} \right)^q
\\ \label{L1final2}
&\le & C\ R^{   (n-2m)(pq-1) -[(a+2m)q+(b+2m)] }  \int_{B_{2R}\setminus B_{R}}      |x|^{b}  u^q      \phi^s_R.
 \end{eqnarray}

\hfill $\Box$

By using the H\"{o}lder's inequality, we can now get the following interpolation estimates on $u$ and $v$.

\begin{cor} \label{L1est}
With the same assumptions as Lemma \ref{L1est},  the following holds. 
\begin{enumerate}
\item[(i)] For any $0<t<q$ and  any $0<s<p$ 
\begin{eqnarray*}
\int_{B_R\setminus B_{R/2}}    v^s \le  C  R^{n-\frac{(a+2m)q+(b+2m)}{pq-1} s} \ \ \text{and} \ \ 
\int_{B_R \setminus B_{R/2} }    u^t \le  C  R^{n-\frac{(b+2m)p+(a+2m)}{pq-1} t}. 
 \end{eqnarray*}
 
 \item[(ii)] For any $0<t<\frac{nq}{n+b}$ and any $0<s<\frac{np}{n+a}$ \begin{eqnarray*}
\int_{B_R}    v^s \le  C  R^{n-\frac{(a+2m)q+(b+2m)}{pq-1} s}\ \ \text{and} \ \ 
\int_{B_R}    u^t \le  C  R^{n-\frac{(b+2m)p+(a+2m)}{pq-1} t}. 
 \end{eqnarray*}
\end{enumerate} 
where $C=C(n,m,t,s,a,b,p,q)$ is independent of $R>1$. 
\end{cor}

We now recall the following fundamental elliptic estimates. We shall apply these estimates frequently for the solutions of (\ref{main}).   

\begin{lemma} (Sobolev inequalities on the sphere $S^{n-1}$) \label{sobolev} 
Let $n\ge2$,  integer $s\ge 1$ and $1<t<\tau\le\infty$. For $z\in W^{s,t}(S^{n-1})$, we have
     $$||z||_{L^\tau(S^{n-1})}\le C || D_\theta^s  z||_{L^t(S^{n-1})}  + C || z  ||_{L^1(S^{n-1})} ,$$ where
     
       $$\left\{
                      \begin{array}{ll}
                        \frac{1}{\tau}= \frac{1}{t}-\frac{s}{n-1}, & \hbox{if $st+1<n$,} \\
                        \tau=\infty, & \hbox{if $st+1>n$,} 
                                                                       \end{array}
                    \right.$$
                     and $C=C(s,t,n,q)>0$  does not depend on $R>1$.
\end{lemma}

\begin{lemma} \label{ellip} 
(Elliptic $L^\tau$-estimate on $B_R$).  Let $1< \tau<\infty$ and $m\ge 1$. For $z\in W^{2m,\tau}(B_{2R})$, we have
$$\int_{B_R\setminus B_{R/2}}   |  D_{x}^{2m} z |^\tau\le C \int_{B_{2R}\setminus B_{R/4}}   |\Delta_x^m z|^\tau  +  C R^{-2m\tau}    \int_{B_{2R}\setminus B_{R/4}}     |z|^\tau ,$$
where $C=C(n,m,\tau)>0$  does not depend on $R>1$.
\end{lemma}

\begin{lemma}\label{interp} 
(An interpolation inequality on $B_R$).  Let $1 \le \tau<\infty$, $m\ge1$ and $1\le i\le 2m-1$.  For sufficiently regular $z$, we have
$$\int_{B_R\setminus B_{R/2}}     |  D_{x}^i z |^\tau \le C  R^{(2m-i)\tau} \int_{B_{2R}\setminus B_{R/4}}   |\Delta_x^{m} z|^\tau  +   C R^{-i\tau}   \int_{B_{2R}\setminus B_{R/4}}    |z|^\tau  ,$$
where $C=C(n,m,\tau,i)>0$  does not depend on $R>1$.
\end{lemma}

By applying Lemma \ref{weightedest}, Corollary \ref{L1est}, Lemma \ref{ellip} and Lemma \ref{interp}, we obtain the following estimates on the derivatives of $u$ and $v$.
\begin{lemma} \label{DL1est}
Let $m\ge 1$ and suppose that either $0\le i\le 2m-1$ and $\epsilon\ge 0$ or $i= 2m$ and  $\epsilon> 0$. Then, for a bounded nonnegative solution $(u,v)$ of (\ref{main}) we have
\begin{eqnarray*}
\int_{B_{2R}\setminus B_{R}}    {    |  D^i_{x}u|^{1+\epsilon}}  &\le& C\  R^{n-i-\frac{(b+2m)p+(a+2m)}{pq-1} +\epsilon(2m-i+a)},
\\
\int_{B_{2R}\setminus B_{R}}      {    |  D^i_{x}v|^{1+\epsilon}}  &\le& C\  R^{n-i-\frac{(a+2m)q+(b+2m)}{pq-1} +\epsilon(2m-i+b)},
 \end{eqnarray*}
where the constant $C=C(a,b,n,m,p,q,i,\epsilon)>0$ does not depend on $R>1$.
\end{lemma}

\noindent \textbf{Proof:}
 For the case  $\epsilon\ge 0$ and $0\le i\le 2m-1$  we apply Lemma \ref{weightedest},  Lemma \ref{interp} and Corollary \ref{L1est} and for the case $\epsilon> 0$ and $i= 2m$ we apply Lemma \ref{weightedest},  Lemma \ref{ellip} and Corollary \ref{L1est}  to get the following 
\begin{eqnarray*}
\int_{B_{2R}\setminus B_{R}}   |D^{i}_{x}u|^{1+\epsilon} &\le& C R^{(2m-i)(1+\epsilon)} \int_{B_{4R}\setminus B_{R/2}}   |\Delta^m u|^{1+\epsilon}  +C\ R^{-i(1+\epsilon)} \int_{B_{4R}\setminus B_{R/2}}   u^{1+\epsilon} \\
&\le& C R^{(2m-i)(1+\epsilon) +  a\epsilon}  \int_{B_{4R}\setminus B_{R/2}}  |x|^{a} v^{p(1+\epsilon)} + C\ R^{-i(1+\epsilon)} \int_{B_{4R}\setminus B_{R/2}}    u^{1+\epsilon} \\
&\le& C R^{(2m-i)(1+\epsilon) +  a\epsilon} \int_{B_{4R}\setminus B_{R/2}}   |x|^{a} v^{p} + C\ R^{-i(1+\epsilon)}  \int_{B_{4R}\setminus B_{R/2}}    u \\ &\le&   C R^{(2m-i)(1+\epsilon) +  a\epsilon} R^{n -2m- \frac{(b+2m)p+(a+2m)}{pq-1}} +C \ R^{-i(1+\epsilon)}  R^{n-\frac{(a+2m)q+(b+2m)}{pq-1}}\\
&\le &C\ R^{n -i- \frac{(b+2m)p+(a+2m)}{pq-1}+\epsilon(2m+a-i)}. 
\end{eqnarray*}
Note that we have used the boundedness assumption on $u$ and $v$ in the above when $\epsilon>0$.  The proof of the other integral estimate on the  gradients of $v$ is quite similar. 

\hfill $\Box$

To prove our main results we apply the following Pohozaev identity.    

\begin{lemma}\label{Poho}
(Pohozaev identity). Suppose $\lambda,\gamma\in\mathbb{R}$ satisfy $\lambda+\gamma=n-2m$.  If $(u,v)$ is a nonnegative solution of (\ref{main}), then it necessarily satisfy
\begin{enumerate}
\item For $m=2k+1$ where $k\ge 0$, 
\begin{eqnarray*}
&&\left(\frac{n+a}{p+1} -\lambda \right) \int_{B_{R}} |x|^{a} v^{p+1} + \left(\frac{n+b}{q+1} -\gamma \right) \int_{B_{R}} |x|^{b} u^{q+1} \\&&
=  \frac{1}{p+1}    \int_{\partial B_R} |x|^a v^{p+1} x\cdot\nu  + \frac{1}{q+1}    \int_{\partial B_R} |x|^b u^{q+1} x\cdot\nu - \int_{\partial B_R} \nabla \Delta ^k u \cdot \nabla \Delta ^k v    x\cdot\nu \\&& + (\lambda+m-1)   \int_{\partial B_R} \Delta ^k v \partial_\nu \Delta^k u + (\gamma+m-1)   \int_{\partial B_R} \Delta ^k u \partial_\nu \Delta^k v
\\&&  + \int_{\partial B_R} \nabla \Delta ^k u\cdot\nu \ x\cdot\nabla \Delta ^k v  + \int_{\partial B_R} \nabla \Delta ^k v\cdot\nu \ x\cdot\nabla \Delta ^k u 
 \\&& + \lambda I(u,v) +  \gamma I(v,u) +   J(u,v)  + J(v,u) 
\end{eqnarray*}
 \item For $m=2k$ where $k\ge 1$,
 \begin{eqnarray*}
&&\left(\frac{n+a}{p+1} -\lambda \right) \int_{B_{R}} |x|^{a} v^{p+1} + \left(\frac{n+b}{q+1} -\gamma \right) \int_{B_{R}} |x|^{b} u^{q+1} \\&&
=  \frac{1}{p+1}    \int_{\partial B_R} |x|^a v^{p+1} x\cdot\nu  + \frac{1}{q+1}    \int_{\partial B_R} |x|^b u^{q+1} x\cdot\nu - \int_{\partial B_R}  \Delta ^k u \Delta ^k v    x\cdot\nu 
 \\&& - \lambda I(u,v) - \gamma I(v,u) -   J(u,v)  - J(v,u) 
\end{eqnarray*}
 \end{enumerate}
where 
\begin{eqnarray*}
 I(u,v)&:=&\sum_{i=0}^{k-1} \int_{\partial B_R} \left( \Delta ^i v \partial_\nu \Delta^{m-i-1} u -\Delta ^{m-i-1} u \partial_\nu \Delta ^i v \right)
 \\ 
 J(u,v)&:=& \sum_{i=0}^{k-1} \int_{\partial B_R}\left(\Delta ^i (x\cdot \nabla v)   \partial_\nu \Delta^{m-i-1} u -\Delta^{m-i-1} u \partial_\nu \Delta ^i (x \cdot\nabla v) \right)
\end{eqnarray*}
\end{lemma} 

\noindent \textbf{Proof:}  The proof is quite standard. We mention few technical facts that facilitates the computations.  Suppose $z,w$ are smooth functions then  for any  $i\in\mathbf N$
\begin{eqnarray*}
\Delta^i(x\cdot \nabla z)&=&2i \Delta^i z+x\cdot \nabla \Delta^i z\\
\nabla z\cdot \nabla (x\cdot\nabla w)+\nabla w\cdot \nabla (x\cdot\nabla z)&=& 2 \nabla z\cdot\nabla w+ x\cdot(\nabla z\cdot\nabla w).
\end{eqnarray*}
Also for any $\lambda,\gamma\in \mathbb R$  the following equalities hold
\begin{eqnarray*}
(\lambda+\gamma) \int_{B_R} \nabla \Delta^ku\cdot\nabla \Delta^k v&=& \lambda \int_{B_R} |x|^a v^{p+1} + \gamma \int_{B_R} |x|^b u^{q+1} + \lambda I(u,v) +\gamma I(v,u) \ \ \text{and} \\
(n-2m) \int_{B_R} \nabla \Delta^ku\cdot\nabla \Delta^k v&=& \int_{B_R} (\Delta ^m u x\cdot \nabla v + \Delta ^m v x\cdot \nabla u) \\&&+ \int_{\partial B_R} \nabla \Delta^ku\cdot\nabla \Delta^k v x\cdot \nu - J(u,v)- J(v,u)
\end{eqnarray*}  
when $m=2k+1$  and for the case $m=2k$ we have similar equations as  
\begin{eqnarray*}
(\lambda+\gamma) \int_{B_R} \Delta^k u  \Delta^k v &=& \lambda \int_{B_R} |x|^a v^{p+1} + \gamma \int_{B_R} |x|^b u^{q+1} - \lambda I(u,v) -\gamma I(v,u) \ \ \text{and} \\
(n-2m) \int_{B_R}  \Delta^ku \Delta^k v&=& -\int_{B_R} (\Delta ^m u x\cdot \nabla v - \Delta ^m v x\cdot \nabla u) \\&&+ \int_{\partial B_R}  \Delta^ku \Delta^k v x\cdot \nu + J(u,v) + J(v,u).
\end{eqnarray*}
Finally for either  $m=2k$ or  $m=2k+1$ we have 
  \begin{eqnarray*}
 \int_{B_R}  |x|^a v^p x\cdot \nabla v +    \int_{B_R}  |x|^b u^q x\cdot \nabla u &=& -\frac{n+a}{p+1} \int_{B_R} |x|^a v^{p+1}  -\frac{n+b}{q+1} \int_{B_R} |x|^b u^{q+1} \\&&+ \frac{1}{p+1} \int_{\partial B_R} |x|^a v^{p+1} x\cdot \nu  +\frac{1}{q+1} \int_{\partial B_R} |x|^b u^{q+1} x\cdot\nu
\end{eqnarray*}
  
\hfill $\Box$

We are now in the position to prove Theorem \ref{result}. The main technique here is to apply the Pohozaev identity that Lemma \ref{Poho} and then taking the advantage of the elliptic regularity theory and in particular the lemmata mentioned before to get certain decay estimates on each boundary term appeared in the Pohozaev identity. 
\\

\noindent \textbf{Proof of Theorem \ref{result}:} Since $(p,q)$ satisfy (\ref{curve}), then we can choose $\lambda$ and $\gamma$ such that $\frac{n+a}{p+1}>\lambda$ and $\frac{n+b}{q+1}>\gamma$. Now, for all $R>1$ define the following positive function of $R$ that is in fact the left hand side of the Pohozaev identity
 $$L(R):=\left(\frac{n+a}{p+1} -\lambda \right) \int_{B_{R}} |x|^{a} v^{p+1} + \left(\frac{n+b}{q+1} -\gamma \right) \int_{B_{R}} |x|^{b} u^{q+1}.$$ 
From Lemma \ref{Poho} and for either $m=2k $ or $m=2k+1$ we have the following upper bound on $L$ 
\begin{equation}\label{F}
L(R)\le C \sum_{i=1}^{3} U_{i}(R),
\end{equation}
where  $C=C(n,m,p,q)$ is independent of $R$ and 
\begin{eqnarray*}
U_{1}(R)&:=& R^{n+a}\int_{S^{n-1}} v^{p+1}(R,\theta) + R^{n+b}\int_{S^{n-1}} u^{q+1}(R,\theta)   
\\
U_{2}(R)&:=&  R^{n-1} \sum_{j=0}^{m-1} \int_{S^{n-1}} \left( |D_x^{j} v(R,\theta)|  |D_x^{2m-j-1} u(R,\theta)|  +  |D_x^{j} u(R,\theta)|  |D_x^{2m-j-1} v(R,\theta)| \right)   \\
U_{3}(R)&:=&  R^{n} \sum_{j=0}^{m-1} \int_{S^{n-1}}  \left(  |D_x^{j+1} v(R,\theta)|  |D_x^{2m-j-1} u(R,\theta)| + |D_x^{j+1} u(R,\theta)|  |D_x^{2m-j-1} v(R,\theta)| \right)
 \end{eqnarray*}
 To get this upper bound we have used the following facts.
 \begin{eqnarray*}
 I(u,v)& \le &k \sum_{i=0}^{k-1} \int_{\partial B_R} \left( |D_x^{2i} v|  |D_x^{2m-2i-1} u| + |D_x ^{2m-2i-2} u| |D_x ^{2i+1} v| \right) \\ & \le &k \sum_{j=0}^{2k-1} \int_{\partial B_R}  |D_x^{j} v|  |D_x^{2m-j-1} u|
\end{eqnarray*}
Note that for any $i\in\mathbb{N}$ we have $\Delta^i(x\cdot \nabla v)=2i \Delta^i v+x\cdot \nabla \Delta^i v$  and also $\nu\cdot \nabla(x\cdot\nabla \Delta^i v)=\nu\cdot\nabla\Delta^i v+  \sum_{s,t}^{n} x_s\nu_t \partial_{x_s x_t} (\Delta^i v) $. Using this we get 
\begin{eqnarray*}
J(u,v)&\le& CI(u,v)+ k \sum_{i=0}^{k-1} \int_{\partial B_R}\left( x\cdot \nabla \Delta^i v \partial_\nu \Delta^{m-i-1} u -\Delta^{m-i-1} u \partial_\nu (x\cdot \nabla \Delta^i v) \right)
\\&\le& CI(u,v)+ k R \sum_{i=0}^{k-1} \int_{\partial B_R}\left( |D_x^{2i+1} v |  |D_x^{2m-2i-1} u|  +  |D_x^{2m-2i-2} u| |D_{x}^{2i+2} v| \right)
\\&\le& CI(u,v)+ k R \sum_{j=0}^{2k-1} \int_{\partial B_R}  |D_x^{j+1} v |  |D_x^{2m-j-1} u|  
 \end{eqnarray*}
   In what follows we find upper bounds on each $U_i$ when $1\le i \le 3$.  Let's first fix  $\epsilon>0$ small enough now and then we pick the appropriate value later.  Also, for the sake of simplicity of notations,  throughout the proof, we use the notation $||w||_t$ to show the $L^t(S^{n-1})$ estimates of $w(R,\theta)$ on the sphere that is $||w||_{L^t{(S^{n-1})}}$ or $\left(\int_{S^{n-1}} w^t(R,\theta)\right)^{1/t}$. Here are the upper bounds.  

 Upper bounds for $U_{1}$. Note that from Lemma \ref{sobolev} we have the Sobolev embedding $W^{2m,1+\epsilon}(S^{n-1})\hookrightarrow L^\infty(S^{n-1})$ in dimension $n=2m+1$. Therefore, 
\begin{eqnarray*}
\left(\int_{S^{n-1}} v^{p+1}(R,\theta)\right)^{\frac{1}{p+1}}&=&|| v||_{p+1}\le || v||_{\infty} \le C || D_{\theta}^{2m} v||_{1+\epsilon}+ C ||v||_{1}  \\&\le&  C R^{2m}|| D_{x}^{2m} v||_{1+\epsilon}+C ||v||_{1}
 \end{eqnarray*}
 So, applying the same argument for $u$ we get
\begin{eqnarray}
\nonumber   U_{1}(R) &\le& C  R^{n+a} \left( R^{2m} || D_{x}^{2m} v||_{1+\epsilon} + || v ||_{1}\right)^{p+1}\\ \label{U1}
&& +C  R^{n+b} \left( R^{2m} || D_{x}^{2m} u||_{1+\epsilon} + || u ||_{1}\right)^{q+1}.
\end{eqnarray}
Upper bounds for $U_{2} $.  For any $j=0,\cdots,m-1$  we have  $1\le j+1 \le m$ and also $1\le 2m-1\le 2m-j-1\le m$.  So  from H\"{o}lder's inequality we get 
\begin{equation*}
\int_{S^{n-1}} |D_x^{j} v(R,\theta)|  |D_x^{2m-j-1} u(R,\theta)| \le  || D^j_x v||_{{\frac{2m}{j+1}}} ||  D_x^{2m-j-1} u||_{{\frac{2m}{2m-j-1}}}
\end{equation*}
Note that  from Lemma \ref{sobolev} we get the embeddings  $W^{2m-j-1,1+\epsilon}(S^{n-1})\hookrightarrow L^{\frac{2m}{j+1}}(S^{n-1})$ and $W^{j+1,1+\epsilon}(S^{n-1})\hookrightarrow L^{\frac{2m}{2m-j-1}}(S^{n-1})$ in dimension $n=2m+1$.  So, 
\begin{equation*}
|| D^j_x v||_{L^{\frac{2m}{j+1}}}   \le C  || D^{2m-j-1}_{\theta}D^{j}_{x}v||_{1+\epsilon} + C ||D^{j}_{x}v||_{1}   \le C  R^{2m-j-1} ||D_{x}^{2m-1}v||_{1+\epsilon} +C ||D^{j}_{x}v||_{1}  
\end{equation*}
and 
\begin{eqnarray*}
\nonumber || D^{2m-j-1}_x u||_{L^{\frac{2m}{2m-j-1}}}   &\le & C   || D^{j+1}_{\theta}D^{2m-j-1}_x u||_{1+\epsilon} +C ||D^{2m-j-1}_x u||_{1}  \\&\le& C R^{j+1} ||D_{x}^{2m} u||_{1+\epsilon} + C ||D^{2m-j-1}_x u||_{1}  
\end{eqnarray*}
Therefore
\begin{eqnarray}\label{U2}
\nonumber  U_{2}(R) &\le& C R^{n-1} \sum_{j=0}^{m-1} \left( R^{j+1} ||D_{x}^{2m} u||_{1+\epsilon} +||D^{2m-j-1}_x u||_{1}  \right) \left( R^{2m-j-1} ||D_{x}^{2m-1}v||_{1+\epsilon} +||D^{j}_{x}v||_{1}  \right) \\&& +  C R^{n-1} \sum_{j=0}^{m-1} \left( R^{j+1} ||D_{x}^{2m} v||_{1+\epsilon} +||D^{2m-j-1}_x v||_{1}  \right) \left( R^{2m-j-1} ||D_{x}^{2m-1}u||_{1+\epsilon} +||D^{j}_{x}u||_{1}  \right) 
\end{eqnarray}
Upper bounds for $U_{3} $. Similar arguments and embedding as for $U_2$ can be used for this term as well.  H\"{o}lder's inequality yields  
\begin{equation*}
 \int_{S^{n-1}}  |D_x^{j+1} v(R,\theta)|  |D_x^{2m-j-1} u(R,\theta)|  \le || D^{j+1}_x v||_{L^{\frac{2m}{j+1}}} ||  D_x^{2m-j-1} u||_{L^{\frac{2m}{2m-j-1}}}
\end{equation*}
Again  from the embeddings  $W^{2m-j-1,1+\epsilon}(S^{n-1})\hookrightarrow L^{\frac{2m}{j+1}}(S^{n-1})$ and $W^{j+1,1+\epsilon}(S^{n-1})\hookrightarrow L^{\frac{2m}{2m-j-1}}(S^{n-1})$  we get 
\begin{eqnarray*}
|| D^{j+1}_x v||_{L^{\frac{2m}{j+1}}}  & \le& C  || D^{2m-j-1}_{\theta}D^{j+1}_{x}v||_{1+\epsilon} +C ||D^{j+1}_{x}v||_{1}  \\& \le & C R^{2m-j-1} ||D_{x}^{2m}v||_{1+\epsilon} + C ||D^{j+1}_{x}v||_{1}  
\end{eqnarray*}
and 
\begin{eqnarray*}
\nonumber || D^{2m-j-1}_x u||_{L^{\frac{2m}{2m-j-1}}}   &\le & C   || D^{j+1}_{\theta}D^{2m-j-1}_x u||_{1+\epsilon} +C ||D^{2m-j-1}_x u||_{1}   \\&\le& C  R^{j+1} ||D_{x}^{2m} u||_{1+\epsilon} +C ||D^{2m-j-1}_x u||_{1}  
\end{eqnarray*}
Therefore
\begin{eqnarray}\label{U3}
\nonumber  U_{3}(R) &\le&  C R^{n} \sum_{j=0}^{m-1} \left( R^{j+1} ||D_{x}^{2m} u||_{1+\epsilon} +||D^{2m-j-1}_x u||_{1}  \right) \left( R^{2m-j-1} ||D_{x}^{2m}v||_{1+\epsilon} +||D^{j+1}_{x}v||_{1}  \right) \\&& +  C R^{n} \sum_{j=0}^{m-1} \left( R^{j+1} ||D_{x}^{2m} v||_{1+\epsilon} +||D^{2m-j-1}_x v||_{1}  \right) \left( R^{2m-j-1} ||D_{x}^{2m}u||_{1+\epsilon} +||D^{j+1}_{x}u||_{1}  \right) 
\end{eqnarray}

Now we are ready to show that the  upper bounds on each $U_i(R)$ converges to zero for an appropriate sequence of $R_l$ when $R_l$ converges to infinity.  To construct such a sequence,  for any $j=0,\cdots,2m-1$ and $i=2m-1,2m$  define the following sets where $M$ is a large constant that will be determined later.   
\begin{eqnarray*}
\label{mv}\Gamma^{(j)}_{1}(R) &:=&\left\{r\ \in(R/2,R); \  ||D_x^j v||_{1}> M R^{-j-\frac{(a+2m)q+(b+2m)}{pq-1}}\right\},\\
\label{mu}\Gamma^{(j)}_{2}(R) &:=&\left\{r\ \in(R/2,R); \  ||D_x^j u||_{1}  > M R^{-j-\frac{(b+2m)p+(a+2m)}{pq-1} }  \right\},\\
\label{mD2v}\Gamma^{(i)}_{3} (R,\epsilon)&:=&\left\{r\ \in(R/2,R); \   ||D^i_{x} v ||^{1+\epsilon}_{1+\epsilon}  > M R^{-i-\frac{(a+2m)q+(b+2m)}{pq-1}+\epsilon(2m+b-i)} \right\},\\
\label{mD2u}\Gamma^{(i)}_{4}(R,\epsilon) &:=&\left\{r\ \in(R/2,R); \   ||D^i_{x} u||^{1+\epsilon}_{1+\epsilon} > M R^{-i-\frac{(b+2m)p+(a+2m)}{pq-1} +\epsilon(2m+a-i) }  \right\}.
\end{eqnarray*}
Note that from Lemma  \ref{DL1est}  for either $\epsilon\ge 0$ and $0\le t\le 2m-1$ or $\epsilon> 0$ and $t= 2m$ we have 
\begin{eqnarray*}
 C &\ge& R^{-n+t+\frac{(b+2m)p+(a+2m)}{pq-1} - \epsilon(2m-t+a)} \int_{B_{R}\setminus B_{R/2}}    {    |  D^t_{x}u|^{1+\epsilon}}  
 \\&=&   R^{-n+t+\frac{(b+2m)p+(a+2m)}{pq-1} - \epsilon(2m-t+a)}   \int_{R/2}^{R} ||   D^t_{x}u ||^{1+\epsilon}_{1+\epsilon} r^{n-1} dr
 \end{eqnarray*} 
So, for $i=2m-1,2m$ we have 
\begin{eqnarray*}
 C &\ge&  R^{-n+i+\frac{(b+2m)p+(a+2m)}{pq-1} - \epsilon(2m-i+a)}   \int_{\Gamma^{(i)}_{4}(R,\epsilon) } ||   D^i_{x}u ||^{1+\epsilon}_{1+\epsilon} r^{n-1} dr 
 \\ &\ge& M  R^{-n+i+\frac{(b+2m)p+(a+2m)}{pq-1} - \epsilon(2m-i+a)}  |\Gamma^{(i)}_{4}(R,\epsilon)| R^{n-1} R^{-i-\frac{(b+2m)p+(a+2m)}{pq-1} +\epsilon(2m+a-i) }
 \\&=& M |\Gamma^{(i)}_{4}(R,\epsilon)|  R^{-1} 
 \end{eqnarray*} 
 that is  $|\Gamma^{(i)}_{4}(R,\epsilon)| \le \frac{C R}{M}$. Similarly one can apply the same argument to show that $|\Gamma^{(i)}_{3}(R,\epsilon)| \le \frac{C R}{M}$,  $|\Gamma^{(j)}_{1}(R)| \le \frac{C R}{M}$ and $|\Gamma^{(j)}_{2}(R)| \le \frac{C R}{M}$. Therefore, we can choose $M$ large enough to make sure that 
 \begin{equation}
\sum_{\tilde i=3}^{4} \sum_{i=2m-1}^{2m} |\Gamma^{(i)}_{\tilde i}(R,\epsilon)| + \sum_{\tilde j=1}^{2} \sum_{j=0}^{2m-1}   |\Gamma^{(j)}_{\tilde j}(R)| \le \frac{(4m+4)CR}{M} \le  \frac{R}{3}
 \end{equation}
 Hence, for each $R\ge 1$, we can find 
\begin{equation} \label{hatr}
  R_l\in (R/2,R)\setminus \left \{\cup_{i=2m-1}^{2m} \cup_{\tilde i=3}^{4}  \Gamma^{(i)}_{\tilde i}(R,\epsilon),\cup_{j=0}^{2m-1} \cup_{\tilde j=1}^{2}  \Gamma^{(j)}_{\tilde j}(R)\right\}  \neq\phi.
\end{equation}
Now we use the sequence $  R_l$ to get a decay estimate on each $U_i(R)$ where $1\le i\le 3$. 

Decay estimate on $U_{1}$. From (\ref{U1}) we get 
\begin{eqnarray*}
U_{1} (R_l) &\le& C  R_l^{n+a} \left( R_l^{2m} R_l^{\left(-2m-\frac{(a+2m)q+(b+2m)}{pq-1}+\epsilon b\right) \frac{1}{1+\epsilon}}  + R_l^{-\frac{(a+2m)q+(b+2m)}{pq-1}}  \right)^{p+1}\\
&& + C  R_l^{n+b} \left( R_l^{2m} R_l^{\left(-2m-\frac{(b+2m)p+(a+2m)}{pq-1}+\epsilon a\right) \frac{1}{1+\epsilon}}  + R_l^{-\frac{(b+2m)p+(a+2m)}{pq-1}}  \right)^{q+1}
\\&\le&  C \left(  R_l^{-f_{1}(\epsilon)} + R_l^{-\tilde f_{1}(\epsilon)}   \right),
\end{eqnarray*}
where 
\begin{eqnarray*}
f_{1}(\epsilon)= (p+1) \left[  \left(2m+\frac{(a+2m)q+(b+2m)}{pq-1}-b\epsilon\right) \frac{1}{1+\epsilon} -2m- \frac{n+a}{p+1}\right],\\
\tilde f_{1}(\epsilon)= (q+1) \left[  \left(2m+\frac{(b+2m)p+(a+2m)}{pq-1}-a\epsilon\right) \frac{1}{1+\epsilon}  -2m- \frac{n+b}{q+1}\right].
\end{eqnarray*}
Decay estimate on $U_{2}$. From (\ref{U2}) we get 
\begin{eqnarray*}
\nonumber  U_{2}(R_l) 
&\le& C 
R_l^{n-1} \sum_{j=0}^{m-1} \left( R_l^{j+1}R_l^{\left(-2m-\frac{(b+2m)p+(a+2m)}{pq-1}+\epsilon a \right)\frac{1}{1+\epsilon}}  +R_l^{-2m+j+1-\frac{(b+2m)p+(a+2m)}{pq-1} }  \right) 
\\&&
\left( R_l^{2m-j-1} R_l^{ \left(-2m+1-\frac{(a+2m)q+(b+2m)}{pq-1} +\epsilon(b+1)\right) \frac{1}{1+\epsilon}}  +R_l^{-j-\frac{(a+2m)q+(b+2m)}{pq-1}}  \right)
 \\&& 
 +  C 
R_l^{n-1} \sum_{j=0}^{m-1} \left( R_l^{j+1}R_l^{\left(-2m-\frac{(a+2m)q+(b+2m)}{pq-1}+\epsilon b \right)\frac{1}{1+\epsilon}}  +R_l^{-2m+j+1-\frac{(a+2m)q+(b+2m)}{pq-1} }  \right) 
\\&&
\left( R_l^{2m-j-1} R_l^{ \left(-2m+1-\frac{(b+2m)p+(a+2m)}{pq-1} +\epsilon(a+1)\right) \frac{1}{1+\epsilon}}  +R_l^{-j-\frac{(b+2m)p+(a+2m)}{pq-1}}  \right)
 \\&=& 2m C  R_l^{n+2m}  \left( R_l^{\left(-2m-\frac{(b+2m)p+(a+2m)}{pq-1}+\epsilon a \right)\frac{1}{1+\epsilon}} +R_l^{-2m-\frac{(b+2m)p+(a+2m)}{pq-1}}  \right)\\&& \left( R_l^{\left(-2m-\frac{(a+2m)q+(b+2m)}{pq-1}+\epsilon b \right)\frac{1}{1+\epsilon}} +R_l^{-2m-\frac{(a+2m)q+(b+2m)}{pq-1}}  \right) 
  \\&\le& C R_l^{-f_2(\epsilon)}
 \end{eqnarray*}
where $f_2(\epsilon)$ is defined as
\begin{equation}\label{f2}
 f_2(\epsilon):= -n +2m \left(\frac{1-\epsilon}{1+\epsilon}\right)+    \frac{(b+2m)(p+1)+(a+2m)(q+1)}{(pq-1)(1+\epsilon)} - \frac{(a+b)\epsilon  }{1+\epsilon}.
 \end{equation}
 Similarly,  from (\ref{U3}) one can show that $$U_3(R_l) \le C R_l^{-f_2(\epsilon)}.$$ From (\ref{F}) and the upper bounds on each $U_i$ we have 
\begin{equation*}
L(R_l)\le C \sum_{i=1}^{3} U_{i}(R_l) \le C \left( R_l^{-f_2(\epsilon)} + R_l^{-f_1(\epsilon)} + R_l^{-\tilde f_1(\epsilon)} \right)
\end{equation*}
For each $\epsilon\ge 0$ define $f(\epsilon):=\min\{f_1(\epsilon),\tilde f_1 (\epsilon),f_2(\epsilon)\}$. So, $L(R_l)\le CR_l^{ -f(\epsilon)}$. Now to finish the proof we show that for $\epsilon > 0$ small enough $f(\epsilon)>0$. Note that by a straightforward calculation, one can see that $\frac{n+a}{p+1}+\frac{n+b}{q+1}>n-2m$ is equivalent to each one of the following inequalities.  
\begin {eqnarray}\label{ineq000} 
f_{2}(0)&=& -n+2m   + \frac{(b+2m)(p+1)+(a+2m)(q+1)}{pq-1}  >0,\\
 \label{ineq111} f_{1}(0)&=& (p+1)  \left(\frac{(a+2m)q+(b+2m)}{pq-1} - \frac{n+a}{p+1} \right)>0 \\
\label{ineq222} \tilde f_{1}(0)&=& (q+1)  \left( \frac{(b+2m)p+(a+2m)}{pq-1} - \frac{n+b}{q+1} \right) >0, 
\end{eqnarray}
Therefore, we can choose $\epsilon>0$ small enough such that $f(\epsilon)>0$. We now conclude by sending $R\to \infty$ that $L(R_l)=0$ and then $u=v=0$. 

\hfill $\Box$

\section{Liouville theorems for radial solutions via ODE arguments}\label{rad}
In this section we focus on the radial solutions of (\ref{main}) and we prove Theorem \ref{resultrad}.   When we are dealing with radial solutions, the weight functions $|x|^a$ and $|x|^b$ would not change the level of difficulty of the problem much. In other words, the methods and ideas that are used for the case $a=b=0$, can be directly adjusted.  Therefore, we omit some of the proofs. What we would like to emphasize in this section is how the  radial assumption make it easier to get decay estimates on solutions of (\ref{main}), see Lemma \ref{hello} and Corollary \ref{coruv}.  Since we do not need to apply Sobolev embeddings and regularity theory, there will be no restriction on the dimension.       The methods that we apply here are strongly motivated by the methods used in \cite{mi,m}. 

\begin{lemma}\label{sign} 
Suppose that $(u,v)$ is a positive solution of (\ref{main}), then $(-\Delta)^{i}u>0$ and $(-\Delta)^{i}v>0$ where $i=1,2,\cdots,m$. 
\end{lemma}

\noindent \textbf{Proof:} The proof directly follows the methods given in \cite{wx} for polyharmonic equations that is also used in \cite{lgz}  for polyharmonic systems.  The idea is to define the average function on $\partial B_R$ as it is defined in \cite{n}. 

\hfill $\Box$

    \begin{lemma}\label{hello} Let $n \ge 3$.  Suppose that $(u,v)$ is a positive radial solution of (\ref{main}). Then the following pointwise decay estimates hold for any $i=0,\cdots,m$ and $j=0,\cdots,m-1$ provided $r>0$
    \begin {eqnarray*}
(-\Delta)^i u(r) &\le& C_{n,m,i} r^{-2i- \frac{(b+2m)p+(a+2m)}{pq-1} } \\
(-\Delta)^i v(r) &\le& C_{n,m,i} r^{-2i- \frac{(a+2m)q+(b+2m)}{pq-1} } \\
| \Delta ^j u' (r) | &\le& C_{n,m,j} r^{-2j-1- \frac{(b+2m)p+(a+2m)}{pq-1}} \\
| \Delta ^j v'(r) | &\le& C_{n,m,j} r^{-2j-1- \frac{(a+2m)q+(b+2m)}{pq-1} }
     \end{eqnarray*}
            \end{lemma}
            
            \noindent \textbf{Proof:}   Define $u_i=(-\Delta)^i u$  where $i=0,\cdots,m-1$. From Lemma \ref{sign} we have $u_i>0$ and $-\Delta u_i > 0$. Therefore, $u_i'<0$.  Note that from the definition of the sequences $(u_i)_i$ we have $u_{i+1}=-\Delta u_i$ when $i=0,\cdots,m-1$ and $r>0$ 
 \begin {eqnarray*} 
-u_i'(r) r^{n-1} &=& \frac{r^n}{n} u_{i+1}(r) -\frac{1}{n} \int_0^r u_{i+1}'(s) s^{n} ds\\
&\ge &  \frac{r^n}{n} u_{i+1}(r)
     \end{eqnarray*}
that is $r u_{i+1}(r) \le - n u_i'(r)$.   On the other hand, since  $u_i>0$ and $-\Delta u_i > 0$  we have $r u_i'+(n-2)u_i \ge 0$ that is $-r u_i' \le (n-2) u_i$.  Therefore, 
\begin{equation}\label{ui+1}
 u_{i+1}(r) \le \left(  \frac{n(n-2)}{r^2} \right)^{i+1}  u \ \ \text{for all } i=0,\cdots,m-1
 \end{equation} 
In particular, $ u_{m}(r) \le \left(  \frac{n(n-2)}{r^2} \right)^{m}  u(r) $  and note that  $u_m(r)=(-\Delta )^m u(r)= r^a v^p(r)$.  Therefore, 
\begin{equation}\label{v}
v^p(r) \le (n(n-2))^m r^{-2m-a} u(r)
\end{equation} 
Similarly, for $v$ we get 
\begin{equation}\label{u}
u^q(r) \le (n(n-2))^m r^{-2m-b} v(r)
\end{equation} 
From (\ref{u}) and  (\ref{v}) we get 
 \begin {eqnarray*} 
u(r) &\le& (n(n-2))^{m\left(\frac{p+1}{pq-1}\right)} r^{-\frac{(a+2m)+(b+2m)p}{pq-1}}  \\
v(r) &\le& (n(n-2))^{m\left(\frac{q+1}{pq-1}\right)} r^{-\frac{(b+2m)+(a+2m)q}{pq-1}} 
     \end{eqnarray*}
Then from (\ref{ui+1}) and (\ref{u})  we have for all $i=0,\cdots,m$
 \begin {eqnarray*} 
u_i=(-\Delta )^i u(r) &\le& (n(n-2))^{i+m\left(\frac{p+1}{pq-1}\right)} r^{-2i-\frac{(a+2m)+(b+2m)p}{pq-1}}  \\
v_i=(-\Delta )^i v(r) &\le& (n(n-2))^{i+m\left(\frac{q+1}{pq-1}\right)} r^{-2i-\frac{(b+2m)+(a+2m)q}{pq-1}} 
     \end{eqnarray*}
To get the other bounds on the derivative of $u_i$ and $v_i$ one can use  $0\le -r u_i' \le (n-2) u_i$ where $i=0,\cdots,m-1$. 

\hfill $\Box$

\begin{cor}\label{coruv}   Let $n \ge 2m+1$ for $m \ge 1$.  Suppose that $(u,v)$ is a positive radial solution of (\ref{main}) and (\ref{curve}) holds. Then for any $t,\tilde t=0,\cdots,m$ and $s,\tilde s=0,\cdots,m-1$ when $R\to\infty$ 
 \begin {eqnarray} 
\label{uv} R^n |\Delta^t u(R)| |\Delta^{\tilde t} v(R)| &\to& 0 \ \ \text{where  }  t+\tilde t=m \\
\label{uv'} R^{n-1} |\Delta^t u(R)| | \Delta ^s v'(R) | &\to& 0 \ \ \text{where  }  t+s=m-1 \\ 
\label{u'v} R^{n-1} |\Delta^t v(R)| | \Delta ^s u'(R) | &\to& 0 \ \ \text{where  }  t+s=m-1 \\ 
\label{u'v'} R^{n}  | \Delta ^s u'(R) | | \Delta ^{\tilde s} v'(R) | &\to& 0 \ \ \text{where  }  s+\tilde s=m-1
     \end{eqnarray}
\end{cor}
  
 \noindent \textbf{Proof of Theorem \ref{resultrad}:}    The idea is to apply the Pohozaev identity as in the proof of Theorem \ref{result}. Since $(p,q)$ satisfy (\ref{curve}), then we can choose $\lambda$ and $\gamma$ such that $\frac{n+a}{p+1}>\lambda$ and $\frac{n+b}{q+1}>\gamma$. Now, for all $R>1$ define the following positive function of $R$ that is the left hand side of the Pohozaev identity
 $$L(R):=\left(\frac{n+a}{p+1} -\lambda \right) \int_{B_{R}} |x|^{a} v^{p+1} + \left(\frac{n+b}{q+1} -\gamma \right) \int_{B_{R}} |x|^{b} u^{q+1}.$$ 
From Lemma \ref{Poho} and for either $m=2k $ or $m=2k+1$ it is straightforward to observe that the following upper bound on $L$ holds
\begin{equation}\label{F}
L(R)\le C \sum_{i=1}^{5} U_i(R)
\end{equation}
where   $C=C(m,n,a,b,p,q)$ is independent of $R$ and 
\begin{eqnarray*}
U_{1}(R)&:=& R^{n+a}  v^{p+1}(R) + R^{n+b} u^{q+1}(R)   
\\
U_{2}(R)&:=&  R^{n-1} \sum_{j=0}^{k}  \left( |\Delta^{j} v(R)|  |  \Delta^{m-j-1} u'(R)|  +   |\Delta^{m-j-1} u(R)| | \Delta^{j} v' (R)|  \right)   \\
U_{3}(R)&:=&  R^{n-1} \sum_{j=0}^{k}  \left( |\Delta^{j} u(R)|  |  \Delta^{m-j-1} v'(R)|  +   |\Delta^{m-j-1} v(R)| | \Delta^{j} u' (R)|  \right)   \\
U_{4}(R)&:=&   R^{n} \sum_{j=0}^{k}   |  \Delta^{j} v' (R)|  |  \Delta^{m-j-1} u'(R)|  + R^{n} \sum_{j=0}^{k-1}  |\Delta^{m-j-1} u(R)| | \Delta^{j+1} v(R)| \\
U_{5}(R)&:=&   R^{n} \sum_{j=0}^{k}   |  \Delta^{j} u' (R)|  | \Delta^{m-j-1} v'(R)|  + R^{n} \sum_{j=0}^{k-1}  |\Delta^{m-j-1} v(R)| | \Delta^{j+1} u(R)| 
 \end{eqnarray*}
 Not that to get this upper bound we have used the following facts. Suppose $w,z$ are radial functions defined on a ball then $\partial_\nu w=\partial_r w$, $\nabla w\cdot\nabla z=w_r z_r$ and $x\cdot\nabla w=r w_r$. Therefore, for any $R>1$
 \begin{eqnarray*}
 \frac{1}{p+1}    \int_{\partial B_R} |x|^a v^{p+1} x\cdot\nu &\le& C   R^{n+a}  v^{p+1}(R)\\
 \frac{1}{q+1}    \int_{\partial B_R} |x|^b u^{q+1} x\cdot\nu &\le& C  R^{n+b} u^{q+1}(R)  \\
  \int_{\partial B_R} \nabla \Delta ^k u \cdot \nabla \Delta ^k v    x\cdot\nu &\le&  C R^{n}    |  \Delta^{k} v' (R)|  |  \Delta^{k} u'(R)| \\
   \int_{\partial B_R} \Delta ^k v \partial_\nu \Delta^k u  &\le& C  R^{n-1}   |\Delta^{k} v(R)|  |  \Delta^{k} u'(R)|  \\
   \int_{\partial B_R} \nabla \Delta ^k u\cdot\nu \ x\cdot\nabla \Delta ^k v  &\le& C R^{n}    |  \Delta^{k} v' (R)|  |  \Delta^{k} u'(R)| \\
   I(u,v)&\le& C   R^{n-1} \sum_{j=0}^{k-1}  \left( |\Delta^{j} v(R)|  | \Delta^{m-j-1} u'(R)|  +   |\Delta^{m-j-1} u(R)| | \Delta^{j} v' (R)|  \right)
 \end{eqnarray*}
 To find an upper bound on $J(u,v)$ we apply the fact that $\Delta^j (x\cdot\nabla v)=2j\Delta^j v +x\cdot\nabla\Delta^j v$ and then 
  \begin{eqnarray*}
J(u,v) & \le&  C_{m} I(u,v) + \sum_{j=0}^{k-1} \int_{\partial B_R}\left( (x\cdot \nabla \Delta^j v)   \partial_\nu \Delta^{m-j-1} u -\Delta^{m-j-1} u \partial_\nu (x \cdot\nabla \Delta^j v) \right)
\\&\le& C_{m} I(u,v) +  R^{n} \sum_{j=0}^{k-1}   |  \Delta^{j} v' (R)|  |  \Delta^{m-j-1} u'(R)|  + CR^{n} \sum_{j=0}^{k-1}  |\Delta^{m-j-1} u(R)| | \Delta^{j+1} v(R)| \\&&+ C R^{n-1} \sum_{j=0}^{k-1}  |\Delta^{m-j-1} u(R)| | \Delta^{j} v'(R)| 
  \end{eqnarray*}
Note that in the above we  have also used the fact that 
  \begin{eqnarray*}
\partial_\nu (x \cdot\nabla \Delta^j v) &=& (R \Delta^j v'(R))' =R\Delta^j v''(R)+\Delta^j v'(R) \\& =& R\Delta^j v''(R)+(n-1)\Delta^j v'(R) - (n-2)\Delta^j v'(R) \\&=& R\Delta ^{j+1}v(R) - (n-2)\Delta^j v'(R) \\ &\le & C R |\Delta ^{j+1}v(R)| + C |\Delta^j v'(R)|
 \end{eqnarray*}
 In what follows we apply Corollary \ref{coruv}  to  show that $L(R)\to 0$ as $R\to \infty$.  Note that $U_1(R) \le C R^n v |\Delta^m u| + C R^n u |\Delta^m v|$. Then applying (\ref{uv}) when  $t=m$ and $\tilde t=0$ and also when  $t=0$ and $\tilde t=m$, we get $U_1(R)\to 0$ as $R\to \infty$. From the decay estimates (\ref{u'v}) and (\ref{uv'}) when $t$ and $s$ are set to be $0\le j\le k$ and $0\le m-k-1\le m-j-1\le m-1$  we get $U_2(R)$ and $U_3(R)\to 0$ as $R\to \infty$. Similarly, from the decay estimates  (\ref{uv}) and  (\ref{u'v'})   we get $U_4(R)$ and $ U_5(R)\to 0$ as $R\to \infty$. Therefore, 
\begin{equation}\label{L0}
L(R)\to 0 \ \ \text{as} \ \  R\to \infty
\end{equation}
On the other hand,  multiplying both equations of (\ref{main}) with $v$ and $u$ we have
 \begin{eqnarray*}
\int_{B_R} |x|^a v^{p+1} &= &\int_{B_R} v  (-\Delta)^m u  \\
\int_{B_R} |x|^b v^{q+1} &= &\int_{B_R} u (-\Delta)^m v  
  \end{eqnarray*}
and 
 \begin{eqnarray*}
\int_{B_R} v \Delta^m u  -    \int_{B_R} u \Delta^m v &=& \sum_{i=0}^{m-1} \int_{\partial B_R} \left( \Delta ^i v \partial_{\nu} \Delta^{m-i-1} u - \Delta ^{m-i-1} u \partial_{\nu} \Delta^{i} v \right)
\\&\le& C   R^{n-1}  \sum_{i=0}^{m-1} | \Delta ^i v|  |\Delta^{m-i-1} u' (R)| +| \Delta ^{m-i-1} u(R)|  |\Delta^{i} v' (R)|
  \end{eqnarray*}
Now applying Corollary \ref{coruv} and in fact the decay estimates (\ref{u'v}) when $t=i$ and $s=m-i-1$ and also (\ref{uv'}) when  $t=m-i-1$ and $s=i$ we get  $\int_{B_R} v \Delta^m u  -    \int_{B_R} u \Delta^m v\to 0$ as $R\to\infty$. Therefore, 
$$\int_{B_R} |x|^a v^{p+1} - \int_{B_R} |x|^b v^{q+1}\to 0 \ \ \text{as} \ \ R\to\infty $$
 From this,  (\ref{L0})  and the fact that $\lambda+\gamma=n-2m$ we get the following as $R\to\infty$
  \begin{eqnarray*}
\left(\frac{n+a}{p+1} + \frac{n+b}{q+1} - (n-2m) \right)  \int_{B_{R}} |x|^{a} v^{p+1} &\to 0& \\
\left(\frac{n+b}{q+1} + \frac{n+a}{p+1} - (n-2m) \right)  \int_{B_{R}} |x|^{b} u^{q+1} &\to 0&
  \end{eqnarray*}
From this and the fact that (\ref{curve}) holds we conclude that $u=v=0$. 
   
   \hfill $\Box$

\end{document}